\documentclass[11pt]{article}
\usepackage{amssymb,amsmath,epsfig}
\usepackage{float}
\allowdisplaybreaks
\def\one{{{{\rm 1} \kern -.19em {\rm l}}}}
\def\C{{{{\rm {\mbox{\small l}}} \kern -.50em {\rm C}}}}
\def\R{{{{\rm l} \kern -.15em {\rm R}}}}
\def\N{{{{\rm l} \kern -.15em {\rm N}}}}
\def\E{{{{\rm l} \kern -.15em {\rm E}}}}
\def\P{{{{\rm l} \kern -.15em {\rm P}}}}
\def\Z{{{{\rm Z} \kern -.35em {\rm Z}}}}
\def\1{{{{\rm 1} \kern -.35em {\rm 1}}}}

\begin{document}
\begin{sloppypar}
\vspace*{0cm}
\begin{center}
{\setlength{\baselineskip}{1.0cm}{ {\Large{\bf
DERIVATIVE OF THE JACOBI POLYNOMIALS WITH RESPECT TO THEIR ORDER AND APPLICATIONS TO INDEFINITE INTEGRATION
\\}} }}
\vspace*{1.0cm}
{\large{\sc{Axel Schulze-Halberg}}}
\end{center}
\noindent \\
Department of Mathematics and Actuarial Science, and Department of Physics, Indiana University Northwest, 3400 Broadway,
Gary IN 46408, USA, E-mail: axgeschu@iu.edu \\ \\
\vspace*{.5cm}
\begin{abstract} 
\noindent
We calculate the first derivative of the Jacobi polynomials with respect to their order in explicit form. This derivative is not an elementary function, but 
contains elementary special cases. As an application, we use our result with a recently devised method for resolving a new class of indefinite integrals 
containing Jacobi polynomials.

\end{abstract}


\noindent \\

\section{Introduction}
Special functions play an important role in mathematics and its applications. They typically appear as solutions of higher-order differential equations, which 
serve as models of processes in the sciences, engineering, economics, among many other fields. As such, there is a vast amount of literature on special functions, 
see for example the comprehensive standard work \cite{olver}. For more application-focused texts, the reader may be referred to \cite{ritelli}, \cite{yang}, 
and references therein. Often special functions depend on numerical parameters that in applications represent certain quantities of the underlying model. 
A famous example is the confluent hypergeometric function that governs the solutions of the quantum-mechanical Schr\"odinger equation for the hydrogen 
atom. Here, the parameters of the confluent hypergeometric function are related to the quantum numbers of the system. It is known that many special functions 
depend on their parameters in a differentiable way. The corresponding derivatives with respect to the parameters are commonly called parameter derivatives. 
Besides being interesting for their mathematical properties, these derivatives occasionally appear in applications. Recent examples include 
a study on statistical covariance between two points \cite{marin}, conditional distributions of a dependent variable in economics \cite{rothe}, 
the confluent algorithm of supersymmetric quantum mechanics \cite{bermudez}, and a Wronskian formula for resolving indefinite integrals \cite{xbatint}. 
Parameter derivatives of special functions are documented in e.g. \cite{olver} \cite{brychkovbuch}, see also the digital resources \cite{dlmf} \cite{wolfram}. 
In contrast to differentiating with respect to the actual argument of a special function, it is often a much more difficult task to determine parameter derivatives. 
During the last decades, consistent progress has been made towards this goal, represented by works on hypergeometric-type functions 
\cite{ancarani} \cite{bujar} \cite{bytev}, Legendre functions \cite{brychleg} \cite{szmyt}, Bessel functions \cite{brychkovbessel} \cite{dunster}, polynomials 
\cite{aktas2} \cite{aktas} \cite{aktas3} \cite{aktas4} \cite{sofotasios}, among many others. In the present note, we will focus on Jacobi polynomials. While the derivatives with respect to 
the parameters have been found and discussed \cite{aktas} \cite{froh}, to the best of our knowledge, this is not true for derivatives with respect to the 
order. By expressing the Jacobi polynomials through the hypergeometric function, the order of the polynomials becomes a real-valued variable, such that 
differentiation is possible. The resulting derivatives are then evaluated at nonnegative integer values of the order. As we will see below, these derivatives 
can be stated in closed form. In contrast to the parameter derivatives of Jacobi polynomials, which are known to also be polynomials \cite{froh}, 
the derivatives with respect to the order are in general non elementary. Particular elementary cases exist, they can contain polynomial, logarithmic, 
or inverse trigonometric terms. Besides determining the aforementioned derivatives, we will show that they can be used to resolve certain indefinite integrals 
involving Jacobi polynomials. This application is based on a recently introduced method \cite{xbatint}. The remainder of this note is organized as follows: 
in section 2 we compute the first derivative of the Jacobi polynomials with respect to their order, and state particular cases. Section 3 is devoted to the application 
of the latter derivative, where we determine a closed-form representation of an indefinite integral class involving Jacobi polynomials. Generalizations of 
our result and a discussion is contained in section 4.

\section{Calculation of the derivative with respect to the order}
We will now compute the first derivative of the Jacobi polynomials with respect to their order. To this end, we first first extend the domain of the order 
to the positive real axis by expressing the Jacobi polynomials through the hypergeometric function. In the next step we calculate the derivative 
of the resulting expression, and evaluate it at nonnegative integer values of the order. It should be stressed that the process of differentiating the 
hypergeometric function is not new in general \cite{ancarani}, but in the present case yields a closed-form result. Let $n$ be a nonnegative integer, and assume that $a$, $b$ are real-valued numbers. We define the Jacobi polynomial of order $n$ with parameters 
$a$ and $b$ as
\begin{eqnarray}
P_n^{(a,b)}(x) &=& \frac{1}{n!} \sum\limits_{k=0}^n \frac{(-n)_k~(a+b+n+1)_k~(a+k+1)_{n-k}}{k!} \left(\frac{1-x}{2}\right)^k. \nonumber
\end{eqnarray}
While here $n$ is a discrete variable, we can extend the definition of the Jacobi polynomials, such that the order is defined on the positive half-axis. 
To this end, we use the following known relation to the hypergeometric function for a positive real number $\mu$:
\begin{eqnarray}
P_\mu^{(a,b)}(x) &=& \frac{\Gamma(a+\mu+1)}{\Gamma(\mu+1)~\Gamma(a+1)}~{}_2F_1\left(-\mu,a+b+\mu+1,a+1,\frac{1-x}{2}\right), \label{jachyp}
\end{eqnarray}
where $\Gamma$ and ${}_2F_1$ denote the Gamma and hypergeometric functions, respectively. Note that the definition (\ref{jachyp}) is valid as long as the parameter $a$ does not take negative integer values. 
We will comment further on this restriction below. Next, we observe that the function $\mu \mapsto P_\mu^{(a,b)}(x)$ is 
continuously differentiable on the positive half-axis, since both the Gamma function and the hypergeometric function are. This statement includes the special case 
$\mu=0$ in the sense of a limit. We can now compute the derivative 
of (\ref{jachyp}) with respect to $\mu$ at nonnegative integer values $n$. We find in a first step
\begin{eqnarray}
\left[\frac{\partial}{\partial \mu}~ P_\mu^{(a,b)}(x)\right]_{\mid \mu=n} &=& \Bigg\{
\left[
\frac{\partial}{\partial \mu}~\frac{\Gamma(a+\mu+1)}{\Gamma(\mu+1)~\Gamma(a+1)}\right] {}_2F_1\left(-\mu,a+b+\mu+1,a+1,\frac{1-x}{2}\right)+ \nonumber \\[1ex]
& & \hspace{-2cm} +~\frac{\Gamma(a+\mu+1)}{\Gamma(\mu+1)~\Gamma(a+1)}~\left[\frac{\partial}{\partial \mu}~{}_2F_1\left(-\mu,a+b+\mu+1,a+1,\frac{1-x}{2}\right) \right]
\Bigg\}_{\mid \mu=n}. \label{jachyp1}
\end{eqnarray}
The derivative in the first line can be found in a straightforward manner. We have
\begin{eqnarray}
\left[\frac{\partial}{\partial \mu}~ \frac{\Gamma(a+\mu+1)}{\Gamma(\mu+1)~\Gamma(a+1)} \right]_{\mu=n}&=& 
\frac{\Gamma(1+a+n)~[\psi(1+a+n)-\psi(1+n)]}{\Gamma(1+a)~\Gamma(1+n)}, \nonumber
\end{eqnarray}
where $\psi$ denotes the digamma function. It remains to determine the derivative of the hypergeometric function in (\ref{jachyp1}), which can be done by 
employing its series representation. We obtain
\begin{eqnarray}
\left[\frac{\partial}{\partial \mu}~{}_2F_1\left(-\mu,a+b+\mu+1,a+1,\frac{1-x}{2}\right) \right]_{\mid \mu=n} &=& \nonumber \\[1ex]
& & \hspace{-9cm} =~
\left[\frac{\partial}{\partial \mu}~\sum\limits_{k=0}^\infty \frac{(-\mu)_k~(a+b+\mu+1)_k}{k!~(a+1)_k}\left(\frac{1-x}{2}\right)^k \right]_{\mu=n} \nonumber \\[1ex]
& & \hspace{-9cm} =~ -\left\{\sum\limits_{k=0}^\infty \frac{(-\mu)_k~(a+b+\mu+1)_k}{k!~(a+1)_k}~\big[\psi(k-\mu)-\psi(-\mu) \big]\left(\frac{1-x}{2}\right)^k \right\}_{\mu=n} +  \label{line1} \\[1ex]
& & \hspace{-9cm} +~ \left\{\sum\limits_{k=0}^\infty \frac{(-\mu)_k~(a+b+\mu+1)_k}{k!~(a+1)_k}~\big[\psi(a+b+\mu+1+k)-\psi(a+b+\mu+1) \big]\left(\frac{1-x}{2}\right)^k 
\right\}_{\mid \mu=n}. \nonumber \\ \label{line2}
\end{eqnarray}
We will now express these series in closed form by using the fact that $\mu$ is taken at a nonnegative integer value $n$. Inspection of 
the series in line (\ref{line2}) shows that only terms pertaining to $1 \leq k \leq n$ contribute. This is so because for $k=0$ the difference of the digamma functions 
vanishes, while for $k \geq n+1$ the term $(-\mu)_k$ equals zero, once we insert $n$. Also, we will rewrite the series in line (\ref{line1}) by splitting it into 
two parts: the first part contains all terms for $k \leq n$, while the second part contains the remaining terms for $k \geq n+1$. The reason for splitting this series is 
that the difference of the digamma functions must be handled differently for the two cases. In addition, we will apply the following relationship for 
the Pochhammer symbol:
\begin{eqnarray}
(-\mu)_k &=& (-\mu)_n~(-\mu+n+1)_{k-n-1}~(-\mu+n). \nonumber
\end{eqnarray}
Upon implementing the above in (\ref{line1}) and (\ref{line2}), we obtain
\begin{eqnarray}
\left[\frac{\partial}{\partial \mu}~{}_2F_1\left(-\mu,a+b+\mu+1,a+1,\frac{1-x}{2}\right) \right]_{\mid \mu=n} &=& \nonumber \\[1ex]
& & \hspace{-9cm} =~
-\left\{\sum\limits_{k=0}^n \frac{(-\mu)_k~(a+b+\mu+1)_k}{k!~(a+1)_k}~\big[\psi(k-\mu)-\psi(-\mu) \big]\left(\frac{1-x}{2}\right)^k \right\}_{\mu=n} + \nonumber \\[1ex]
& & \hspace{-9cm} +~ \left\{\sum\limits_{k=1}^n \frac{(-\mu)_k~(a+b+\mu+1)_k}{k!~(a+1)_k}~\big[\psi(a+b+\mu+1+k)-\psi(a+b+\mu+1) \big]\left(\frac{1-x}{2}\right)^k
\right\}_{\mid \mu=n} -
\nonumber \\[1ex]
& & \hspace{-9cm} 
-~\Bigg\{\sum\limits_{k=n+1}^\infty \frac{(-\mu)_n~(-\mu+n+1)_{k-n-1}~(a+b+\mu+1)_k}{k!~(a+1)_k}~(-\mu+n)~\big[\psi(k-\mu)-\psi(-\mu) \big] \times \nonumber \\[1ex]
& & \hspace{-9cm} 
~\times \left(\frac{1-x}{2}\right)^k \Bigg\}_{\mid \mu=n}. \nonumber
\end{eqnarray}
We proceed by inserting the value $n$ in place of $\mu$, except in the differences of the digamma functions in the first and third sum. The result reads
\begin{eqnarray}
\left[\frac{\partial}{\partial \mu}~{}_2F_1\left(-\mu,a+b+\mu+1,a+1,\frac{1-x}{2}\right) \right]_{\mid \mu=n} &=& \nonumber \\[1ex]
& & \hspace{-9cm} =~
-\sum\limits_{k=0}^n \frac{(-n)_k~(a+b+n+1)_k}{k!~(a+1)_k}~\lim\limits_{\mu \rightarrow n} \big[\psi(k-\mu)-\psi(-\mu) \big]\left(\frac{1-x}{2}\right)^k + \label{zeile1} \\[1ex]
& & \hspace{-9cm} +~ \sum\limits_{k=1}^n \frac{(-n)_k~(a+b+n+1)_k}{k!~(a+1)_k}~\big[\psi(a+b+n+1+k)-\psi(a+b+n+1) \big]\left(\frac{1-x}{2}\right)^k -
\label{zeile2} \\[1ex]
& & \hspace{-9cm} 
-\sum\limits_{k=n+1}^\infty \frac{(-n)_n~(1)_{k-n-1}~(a+b+n+1)_k}{k!~(a+1)_k}~\lim\limits_{\mu \rightarrow n} (-\mu+n)~\big[\psi(k-\mu)-\psi(-\mu) \big] \times \label{zeile3} \\[1ex]
& & \hspace{-9cm} 
~\times \left(\frac{1-x}{2}\right)^k. \nonumber
\end{eqnarray}
The two limits in the latter function can be evaluated by means of identities for the digamma function. The first of these identities is the difference equation
\begin{eqnarray}
\lim\limits_{\mu \rightarrow n} \big[\psi(k-\mu)-\psi(-\mu)\big] &=& \sum\limits_{j=0}^{k-1} \frac{1}{j-n}, \label{diff}
\end{eqnarray}
that we will apply to the sums in line (\ref{zeile1}) and line (\ref{zeile2}). Note that the lowest-order term in (\ref{diff}) does not contribute, 
and that the sum is defined due to $k \leq n$. The second identity is the following limit, valid for $k > n$:
\begin{eqnarray}
\lim\limits_{\mu \rightarrow n} (-\mu+n)~\big[\psi(k-\mu)-\psi(-\mu) \big]  &=& 1. \label{limi}
\end{eqnarray}
This limit will be used for the sum in line (\ref{zeile3}). Upon implementing (\ref{diff}) and (\ref{limi}) in (\ref{zeile1})-(\ref{zeile3}), we obtain
\begin{eqnarray}
\left[\frac{\partial}{\partial \mu}~{}_2F_1\left(-\mu,a+b+\mu+1,a+1,\frac{1-x}{2}\right) \right]_{\mid \mu=n} &=& \nonumber \\[1ex]
& & \hspace{-5cm} =~
-\sum\limits_{k=1}^n \sum\limits_{j=0}^{k-1} \frac{(-n)_k~(a+b+n+1)_k}{(j-n)~k!~(a+1)_k} \left(\frac{1-x}{2}\right)^k + \label{summe1} \\[1ex]
& & \hspace{-5cm}  +~\sum\limits_{k=1}^n \sum\limits_{j=0}^{k-1} \frac{(-n)_k~(a+b+n+1)_k}{(a+b+n+1+j)~k!~(a+1)_k} \left(\frac{1-x}{2}\right)^k -
\label{summe2} \\[1ex]
& & \hspace{-5cm} 
-\sum\limits_{k=n+1}^\infty \frac{(-n)_n~(1)_{k-n-1}~(a+b+n+1)_k}{k!~(a+1)_k} \left(\frac{1-x}{2}\right)^k. \label{summe3}
\end{eqnarray}
In order to find a closed-form representation for the series in the last line, we identify it with a hypergeometric function. After factoring out terms, 
we arrive at
\begin{eqnarray}
\sum\limits_{k=n+1}^\infty \frac{(-n)_n~(1)_{k-n-1}~(a+b+n+1)_k}{k!~(a+1)_k} \left(\frac{1-x}{2}\right)^k &=& \nonumber \\[1ex]
& & \hspace{-9cm} =~
\frac{(-n)_n~(a+b+n+1)_{n+1}}{2^{n+1}~(n+1)!~(a+1)_{n+1}}~(1-x)^{n+1}~ {}_3F_2\left(
1,1,2+a+b+2~n,n+2,n+a+2,\frac{1-x}{2}\right). \nonumber
\end{eqnarray}
We can now replace the sum in (\ref{summe3}) with the latter expression, yielding the derivative of the hypergeometric function. After joining the 
sums in lines (\ref{summe1}) and (\ref{summe2}), we arrive at
\begin{eqnarray}
\left[\frac{\partial}{\partial \mu}~{}_2F_1\left(-\mu,a+b+\mu+1,a+1,\frac{1-x}{2}\right) \right]_{\mid \mu=n} &=& \nonumber \\[1ex]
& & \hspace{-9cm} =~
\sum\limits_{k=1}^n \sum\limits_{j=0}^{k-1} \frac{(-n)_k~(a+b+n+1)_k}{k!~(a+1)_k} \left(-\frac{1}{j-n}+\frac{1}{a+b+n+1+j} \right)
\left(\frac{1-x}{2}\right)^k  -
\nonumber \\[1ex]
& & \hspace{-9cm} 
-~\frac{(-n)_n~(a+b+n+1)_{n+1}}{2^{n+1}~(n+1)!~(a+1)_{n+1}}~(1-x)^{n+1}~ {}_3F_2\left(
1,1,2+a+b+2~n,n+2,n+a+2,\frac{1-x}{2}\right). \nonumber
\end{eqnarray}
In the last step we plug this derivative into (\ref{jachyp1}). As a result, we find the derivative of the Jacobi polynomial with respect to the order in its 
final form. It reads
\begin{eqnarray}
\left[\frac{\partial}{\partial \mu}~ P_\mu^{(a,b)}(x)\right]_{\mid \mu=n} &=& \nonumber \\[1ex]
& & \hspace{-4cm} =~
\frac{\Gamma(1+a+n)~[\psi(1+a+n)-\psi(1+n)]}{\Gamma(1+a)~\Gamma(1+n)}~{}_2F_1\left(-n,a+b+n+1,a+1,\frac{1-x}{2}\right)+
\nonumber \\[1ex]
& & \hspace{-4cm} +~
\frac{\Gamma(a+n+1)}{\Gamma(n+1)~\Gamma(a+1)} \times \nonumber \\[1ex]
& & \hspace{-4cm} \times~
\Bigg\{\sum\limits_{k=1}^n \sum\limits_{j=0}^{k-1} \frac{(-n)_k~(a+b+n+1)_k}{k!~(a+1)_k} \left(-\frac{1}{j-n}+\frac{1}{a+b+n+1+j} \right)
\left(\frac{1-x}{2}\right)^k  -\nonumber \\[1ex]
& & \hspace{-4cm} -~
\frac{(-n)_n~(a+b+n+1)_{n+1}}{2^{n+1}~(n+1)!~(a+1)_{n+1}}~(1-x)^{n+1}~ {}_3F_2\left(
1,1,2+a+b+2~n,n+2,n+a+2,\frac{1-x}{2}\right)
\Bigg\}. \nonumber \\ \label{jacder}
\end{eqnarray}
We observe that this derivative is in general not elementary due to the occurrence of the hypergeometric function ${}_3F_2$. As mentioned in 
the introduction, this is in contrast to the parameter derivatives of the Jacobi polynomials that are known to have polynomial form. Elementary cases 
of (\ref{jacder}) exist, 
for example if the third argument of the latter function attains certain integer values. The resulting expressions can be polynomial, or involve digamma, 
logarithmic, or inverse trigonometric terms. In order to illustrate this fact, let us list three special cases of the derivative (\ref{jacder}), where for the 
sake of brevity we restrict ourselves to the case $n=0$. Let us first generate a case that can render (\ref{jacder}) in polynomial form, obtained by the 
parameter setting $a=-b-3$. Inspection of (\ref{jacder}) shows that the third argument of the ${}_3F_2$ function becomes $-1$. The derivative itself 
reads
\begin{eqnarray}
\left[\frac{\partial}{\partial \mu}~ P_\mu^{(-b-3,b)}(x)\right]_{\mid \mu=0} &=& \frac{(x-1)~(5+4~b-x)}{4~(b+1)~(b+2)}+\gamma+\psi(-b-2). \nonumber
\end{eqnarray}
A second special case of (\ref{jacder}) is found by setting $a=b=1$, leaving the third argument of the ${}_3F_2$ function equal to 4. We obtain the 
derivative (\ref{jacder}) in the form
\begin{eqnarray}
\left[\frac{\partial}{\partial \mu}~ P_\mu^{(1,1)}(x)\right]_{\mid \mu=0} &=& 
\frac{1+3~x}{2~(x+1)}+\log\left(\frac{1-x}{2} \right). \nonumber
\end{eqnarray}
In the last special case we generate an inverse trigonometric term in (\ref{jacder}) by setting $a=b=-1/2$, such that the third argument of the ${}_3F_2$ function 
becomes $1$. This gives 
\begin{eqnarray}
\left[\frac{\partial}{\partial \mu}~ P_\mu^{(-\frac{1}{2},\frac{1}{2})}(x)\right]_{\mid \mu=0} &=& 
\frac{2~(x-1)}{\sqrt{1-x^2}}~\arcsin\left(\sqrt{\frac{1-x}{2}} \right) -\log(4). \nonumber
\end{eqnarray}
Now that we have obtained our main result (\ref{jacder}), we proceed by employing it to resolve a class of indefinite integrals containing Jacobi polynomials.

\section{Application: Wronskian integral formula}
As mentioned in the introduction, recently a method was devised for generating indefinite integrals through parameter derivatives \cite{xbatint}. For the 
sake of completeness, we will now summarize the main result from the latter reference: suppose that $u$ and $V$ are smooth functions, defined on an open subset of $\mathbb{R}^2$, such that the equation
\begin{eqnarray}
\frac{\partial^2 u(x,\mu)}{\partial x^2}+V(x,\mu)~u(x,\mu) &=& 0, \label{sse}
\end{eqnarray}
is satisfied. Then the following formula for indefinite integration holds:
\begin{eqnarray}
\int\limits^x \frac{\partial V(t,\mu)}{\partial \mu}~ u(t,\mu)^2~dt &=& 
W_{\frac{\partial u}{\partial \mu},u}(x,\mu), \label{int1} 
\end{eqnarray}
where $W$ denotes the Wronskian of the functions in its index. We observe that (\ref{int1}) allows for resolving the indefinite integral on the left side 
by evaluation of the Wronskian on the right side, which merely involves differentiation. Our goal in this application is to match the general form of equation (\ref{sse}) with 
a particular case solved by Jacobi polynomials, where the parameter $\mu$ is given by their order.  As a consequence, the Wronskian on the right side 
of (\ref{int1}) contains the derivative of the Jacobi polynomials with respect to the order $\mu$, which is precisely our result (\ref{jacder}). Let us now start out by 
stating the equation solved by Jacobi polynomials. We recall that it can be written in the form
\begin{eqnarray}
(1-x^2)~\frac{d^2}{dx^2}~P_n^{(a,b)}(x)+[b-a-(a+b+2)~x]~ \frac{d}{dx}~P_n^{(a,b)}(x)+n~(n+a+b+1)~P_n^{(a,b)}(x) = 0. \nonumber
\end{eqnarray}
This equation does not have the form (\ref{sse}), since it contains a first-order derivative term. Therefore, in the next step we will gauge away that term, and we will 
generalize the order $n$ of the Jacobi polynomials to the real variable $\mu$. We define
\begin{eqnarray}
u(x,\mu) &=& (1-x)^\frac{a+1}{2} (1+x)^\frac{b+1}{2}~P_\mu^{(a,b)}(x), \label{u}
\end{eqnarray}
recall that for non-integer values of $\mu$ the definition (\ref{jachyp}) applies. Substitution shows that (\ref{u}) is a solution to
\begin{eqnarray}
& & \frac{\partial^2 u(x,\mu)}{\partial x^2}+\Bigg[
\frac{1-a^2}{4~(x-1)^2}+\frac{1-b^2}{4~(1+x)^2}-\frac{(1+a)~(1+b)-2~\mu~(1+a+b)-2~\mu^2}{4~(x-1)}+ \nonumber \\[1ex]
& & \hspace{5cm}
+~\frac{(a+1)~(b+1)+2~\mu~(1+a+b)+2~\mu^2}{4~(x+1)}
\Bigg]~u(x,\mu) ~=~ 0, \nonumber 
\end{eqnarray}
which matches the general form (\ref{sse}), provided we define the function $V$ as
\begin{eqnarray}
V(x,\mu) &=& 
\frac{1-a^2}{4~(x-1)^2}+\frac{1-b^2}{4~(1+x)^2}-\frac{(1+a)~(1+b)-2~\mu~(1+a+b)-2~\mu^2}{4~(x-1)}+ \nonumber \\[1ex]
&+& \frac{(a+1)~(b+1)+2~\mu~(1+a+b)+2~\mu^2}{4~(x+1)}. \nonumber
\end{eqnarray}
We are now ready to evaluate the left side of the integral formula (\ref{int1}). To this end, we need the derivative of $V$ with respect to $\mu$, which 
can be readily computed as
\begin{eqnarray}
\frac{\partial V(x,\mu)}{\partial \mu} &=& \frac{1+a+b+2~\mu}{1-x^2}. \nonumber
\end{eqnarray}
Substitution of this derivative, along with (\ref{u}), renders the left side of (\ref{int1}) in the following form
\begin{eqnarray}
\int\limits^x \frac{\partial V(t,\mu)}{\partial \mu}_{\mid \mu=n}~ u(t,n)^2~dt = 
(1+a+b+2~n) \int\limits^x (1-t)^a~(1+t)^b~P_n^{(a,b)}(t)^2~dt. \label{inte}
\end{eqnarray}
Observe that we are taking the integrand at a nonnegative integer value $n$ of $\mu$, such that it contains Jacobi polynomials. To the best of our knowledge, 
a closed form of the integral (\ref{inte}) has not been reported in the literature so far. Now that we completed 
evaluating the left side of (\ref{int1}), we continue with the Wronskian on its right side. In explicit form, it reads
\begin{eqnarray}
W_{\frac{\partial u}{\partial \mu},u}(x,n) 
&=& \left[\frac{\partial u(x,\mu)}{\partial \mu}\right]_{\mid \mu=n} \frac{\partial u(x,n)}{\partial x}-
u(x,n) \left[\frac{\partial^2 u(x,\mu)}{\partial x ~\partial \mu} \right]_{\mu=n}. \label{wronskian}
\end{eqnarray}
Let us consider each of the four terms on the right side separately. According to (\ref{u}), the first term can be written as
\begin{eqnarray}
\left[\frac{\partial u(x,\mu)}{\partial \mu}\right]_{\mu=n} &=& (1-x)^\frac{a+1}{2} (1+x)^\frac{b+1}{2} \left[\frac{\partial}{\partial \mu}~P_\mu^{(a,b)}(x) \right]_{\mu=n}. 
\label{dum}
\end{eqnarray}
Note that the derivative on the right side is precisely our result (\ref{jacder}). Hence, upon substituting it into (\ref{dum}), we have the first term contained in 
the Wronskian. The second term can be determined by differentiation in the variable $x$. We find
\begin{eqnarray}
\frac{\partial u(x,\mu)}{\partial x} &=& \frac{1}{2}~(1-x)^\frac{a-1}{2}~(1+x)^\frac{b-1}{2}~\Bigg\{
\Bigg[(1+a+b+\mu)~(x^2-1)~P_{\mu-1}^{(a+1,b+1)}(x)\Bigg]- \nonumber \\[1ex]
&-& \Bigg[a-b+(2+a+b)~x \Bigg]~P_\mu^{(a,b)}(x) \Bigg\}. \label{dux}
\end{eqnarray}
The third term in the Wronskian (\ref{wronskian}) is the function $u$ itself, given in (\ref{u}), where $\mu$ must be replaced by $n$. 
It now remains to compute the mixed derivative in the 
last term of (\ref{wronskian}), which is done by differentiating (\ref{dum}) with respect to $x$. In summary, the Wronskian (\ref{wronskian}) can be 
represented in explicit form by combining the results (\ref{dum}), (\ref{dux}), (\ref{u}), and the derivative of (\ref{dum}) with respect to $x$. While the 
process is straightforward, the resulting expression is very large and involved, such that we omit to state it here. Instead, we will resort to a particular case that 
makes the calculations fairly simple and transparent.

\paragraph{An explicit example.} Let us set $a=b=1$, and $n=2$ in (\ref{inte}) and (\ref{wronskian}). The integral then takes the form
\begin{eqnarray}
7 \int\limits^x (t^2-1)~P_2^{(1,1)}(t)^2~dt &=& \frac{1}{16} \int\limits^x \left(63-693~t^2+2205~t^4-1575~t^6\right) dt.
\nonumber
\end{eqnarray}
While it is clear that this integral can be resolved immediately, we use it here with our formula (\ref{int1}) to demonstrate how it works. 
We plug our settings $a=b=1$, $n=2$ into (\ref{wronskian}), and evaluate the resulting expression. We find after simplification
\begin{eqnarray}
W_{\frac{\partial u}{\partial \mu},u}(x,2) &=& \frac{1}{2}~(x^2-1)^2\Bigg\{-\Bigg[5~P_2^{(1,1)}(x)~\left(\frac{\partial}{\partial \mu}~
P_\mu^{(2,2)}(x) \right)_{\mid \mu=1}+ \nonumber \\[1ex]
&+&
P_{1}^{(2,2)}(x)~\bigg[-P_2^{(1,1)}(x)+5 \left[\frac{\partial}{\partial \mu}~P_\mu^{(1,1)}(x) \right]_{\mid \mu=2} \bigg] \Bigg] \Bigg\}. \label{wronskian1}
\end{eqnarray}
The two derivatives on the right side are special cases of our result (\ref{jacder}). A straightforward computation gives
\begin{eqnarray}
\left(\frac{\partial}{\partial \mu}~P_\mu^{(2,2)}(x) \right)_{\mid \mu=1} &=& \frac{-32+x+134~x^2+97~x^3}{20~(1+x)^2}+3~x~\log\left(\frac{1+x}{2}\right) \nonumber \\[1ex]
\left(\frac{\partial}{\partial \mu}~P_\mu^{(1,1)}(x) \right)_{\mid \mu=2} &=& \frac{-7-59~x+19~x^2+79~x^3}{16~(1+x)}+\frac{3}{4}~(5~x^2-1)~\log\left(\frac{1+x}{2}\right). 
\nonumber 
\end{eqnarray}
We now insert these derivatives into the Wronskian (\ref{wronskian1}). After collecting terms, we arrive at the final form
\begin{eqnarray}
W_{\frac{\partial u}{\partial \mu},u}(x,2) &=& -\frac{3}{16}~\left(16-21~x+77~x^3-147~x^5+75~x^7\right). \label{wronskian2}
\end{eqnarray}
Up to an additive constant, this is the indefinite integral in (\ref{inte}). We can verify directly by either integrating the latter, or by taking the derivative 
of (\ref{wronskian2}).

\section{Discussion and concluding remarks}
Let us now conclude this work with a few comments on the relation (\ref{jachyp}), and on the extension of our main result (\ref{jacder}) to polynomials other 
than the Jacobian case. First, we recall that the identity (\ref{jachyp}) between the Jacobi polynomials and the hypergeometric function is undefined for 
negative integer values of the parameter $a$. This issue can be addressed by employing an alternative representation to (\ref{jachyp}), given by
\begin{eqnarray}
P_\mu^{(a,b)}(x) &=& \frac{\Gamma(-b)}{\Gamma(\mu+1)~\Gamma(-\mu-b)}~{}_2F_1\left(-\mu,a+b+\mu+1,b+1,\frac{1+x}{2}\right), \label{jachyp2}
\end{eqnarray}
which is valid as long as the parameter $b$ does not take integer values. Besides (\ref{jachyp}) and (\ref{jachyp2}), there are further representations 
of the Jacobi polynomials, each of which involve different constraints on the parameters. Next, let us point out that for suitable parameter settings, 
the Jacobi polynomials reduce to a variety of special cases, representing other polynomial systems. As a consequence, our result on the derivative 
(\ref{jacder}) can be used for the latter special cases. For example, we recall the following definition of the Gegenbauer polynomials
\begin{eqnarray}
C_n^a(x) &=& \frac{(2~a)_n}{\left(a+\frac{1}{2} \right)_n}~P_n^{\left(a-\frac{1}{2},a-\frac{1}{2}\right)}(x). \nonumber
\end{eqnarray}
We can calculate the derivative of these polynomials with respect to their order $n$. This gives
\begin{eqnarray}
\left[\frac{\partial}{\partial \mu}~C_\mu^a(x) \right]_{\mid \mu=n} &=& \frac{(2~a)_n~\left[\psi(2~a+n)-\psi\left(a+n+\frac{1}{2} \right) \right]}
{\left(a+\frac{1}{2} \right)_n}~P_n^{\left(a-\frac{1}{2},a-\frac{1}{2}\right)}(x)+ \nonumber \\[1ex]
&+& \frac{(2~a)_n}{\left(a+\frac{1}{2} \right)_n}~
\left[\frac{\partial}{\partial \mu}~P_\mu^{\left(a-\frac{1}{2},a-\frac{1}{2}\right)}(x) \right]_{\mid \mu=n}, \label{gegender}
\end{eqnarray}
where the derivative of the Jacobi polynomials can be found in (\ref{jacder}). We omit to insert the latter result, as (\ref{gegender}) would become 
excessively long. Next, we will compute the derivative of the Legendre polynomials with respect to their order. To this end, we use the following 
relation between them and the Jacobi polynomials:
\begin{eqnarray}
P_n(x) &=& P_n^{(0,0)}(x). \nonumber
\end{eqnarray}
The sought derivative is given by
\begin{eqnarray}
\left[\frac{\partial}{\partial \mu}~P_\mu(x) \right]_{\mid \mu=n} &=& 
\left[\frac{\partial}{\partial \mu}~P_\mu^{(0,0)}(x) \right]_{\mid \mu=n}, \nonumber
\end{eqnarray}
which can be evaluated directly from (\ref{jacder}). As a further example, we will now consider the Chebyshev polynomials of the first kind. These can be 
written in the form
\begin{eqnarray}
T_n(x) &=& \frac{P_n^{\left(-\frac{1}{2},-\frac{1}{2}\right)}(x)}{P_n^{\left(-\frac{1}{2},-\frac{1}{2}\right)}(1)}. \nonumber
\end{eqnarray}
We calculate the derivative of these polynomials with respect to their order as follows:
\begin{eqnarray}
\left[\frac{\partial}{\partial \mu}~T_\mu(x) \right]_{\mid \mu=n} &=& \frac{1}{P_n^{\left(-\frac{1}{2},-\frac{1}{2}\right)}(x)}
\left[\frac{\partial}{\partial \mu}~P_\mu^{\left(-\frac{1}{2},-\frac{1}{2}\right)}(x) \right]_{\mid \mu=n}- \nonumber \\[1ex] &-&
\frac{P_n^{\left(-\frac{1}{2},-\frac{1}{2}\right)}(x)}{P_n^{\left(-\frac{1}{2},-\frac{1}{2}\right)}(1)^2} 
\left[\frac{\partial}{\partial \mu}~P_\mu^{\left(-\frac{1}{2},-\frac{1}{2}\right)}(1) \right]_{\mid \mu=n}. \nonumber
\end{eqnarray}
Again, we observe that the remaining derivatives can be found from our result (\ref{jacder}). 
As a last example, we consider the Chebyshev polynomials of the second kind, given by the definition
\begin{eqnarray}
U_n(x) &=& (n+1)~\frac{P_n^{\left(\frac{1}{2},\frac{1}{2}\right)}(x)}{P_n^{\left(\frac{1}{2},\frac{1}{2}\right)}(1)}. \nonumber
\end{eqnarray}
The derivative of these polynomials with respect to their order is found in a similar way as for their counterparts above. We compute
\begin{eqnarray}
\left[\frac{\partial}{\partial \mu}~U_\mu(x) \right]_{\mid \mu=n} &=& \frac{P_n^{\left(\frac{1}{2},\frac{1}{2}\right)}(x)}{P_n^{\left(\frac{1}{2},\frac{1}{2}\right)}(1)}
+\frac{n+1}{P_n^{\left(\frac{1}{2},\frac{1}{2}\right)}(x)}
\left[\frac{\partial}{\partial \mu}~P_\mu^{\left(\frac{1}{2},\frac{1}{2}\right)}(x) \right]_{\mid \mu=n}- \nonumber \\[1ex]
&-&
\frac{(n+1)~P_n^{\left(\frac{1}{2},\frac{1}{2}\right)}(x)}{P_n^{\left(\frac{1}{2},\frac{1}{2}\right)}(1)^2} 
\left[\frac{\partial}{\partial \mu}~P_\mu^{\left(\frac{1}{2},\frac{1}{2}\right)}(1) \right]_{\mid \mu=n}. \nonumber
\end{eqnarray}
As a final comment, we mention that our integral formula (\ref{int1}) is applicable to the special cases of the Jacobi polynomials mentioned in this section. 
Consequently, the formula can be used to set up and resolve indefinite integrals involving those polynomials.

\paragraph{Disclosure statement.} No potential conflict of interest was reported by the author.

\end{sloppypar}

\end{document}